\documentclass[10pt]{amsart}
\usepackage[style=alphabetic]{biblatex} 
\usepackage{tikz}
\usetikzlibrary{angles}
\usepackage{amssymb}
\usepackage{amsthm}
\usepackage[all]{xy}
\usepackage{microtype}
\usepackage{xcolor}
\usepackage{hyperref}
\usepackage[nameinlink]{cleveref}
\usepackage{graphicx}

\makeatletter
    
    \@addtoreset{equation}{section}
  \makeatother

\hypersetup{
 colorlinks,
 linkcolor={teal},
 citecolor={teal},
 urlcolor={teal}
}

\calclayout

\DeclareFieldFormat
  [article,book,inbook,incollection,inproceedings,patent,thesis,unpublished]
  {title}{\emph{#1\isdot}}

\theoremstyle{plain}
	\newtheorem{theorem}{Theorem}

\theoremstyle{definition}

\begin{document}
\title{Isosceles trapezoids of unit area with vertices in sets of infinite planar measure}
\author[J. Koizumi]{Junnosuke Koizumi}
\address{RIKEN iTHEMS, Wako, Saitama 351-0198, Japan}
\email{junnosuke.koizumi@riken.jp}

\date{\today}
\thanks{}

\begin{abstract}
Paul Erd\H{o}s posed the question of whether every measurable planar set of infinite Lebesgue measure contains the four vertices of an isosceles trapezoid of unit area.
In this paper, we provide an affirmative answer to this question.
Additionally, we present affirmative solutions to similar questions by Erd\H{o}s concerning isosceles triangles and right-angled triangles.
\end{abstract}

\maketitle
\setcounter{tocdepth}{1}
\tableofcontents

\enlargethispage*{20pt}
\thispagestyle{empty}

\section{Introduction}
Paul Erd\H{o}s observed that every measurable subset of the plane of infinite Lebesgue measure contains the vertices of a triangle of area $1$. In the proceedings \cite{Erd83} of the conference \emph{Measure Theory} held in 1983, he posed a number of questions of this type by replacing triangles with the following types of polygons:
\begin{enumerate}
    \item isosceles trapezoids.
    \item cyclic quadrilaterals.
    \item isosceles triangles.
    \item right-angled triangles.
    \item convex polygons with congruent sides.
\end{enumerate}

Recently, Kova\v{c} and Predojevi\'c (\cite{KP}) provided a positive answer to (2) and a negative answer to (5).
However, the questions regarding (1), (3), and (4) appear to remain unresolved to date; see \cite[Problem \#353]{EP}.
In this paper, we provide affirmative answers to all three remaining questions.
In fact, for (3) and (4), we present the following slightly stronger result:

\begin{theorem}\label{thm:iso_right}
    Let $S\subset\mathbb{R}^2$ be an unbounded measurable set of positive Lebesgue measure. Then $S$ contains the vertices of an isosceles triangle of area $1$.
    Also, $S$ contains the vertices of a right-angled triangle of area $1$.
\end{theorem}

Our proof of \Cref{thm:iso_right} is based on the following idea, inspired by the proof in \cite[Theorem 1]{KP}.
Suppose that $O=(0,0)\in S$.
If we define
$$
    f(r\cos\theta,r\sin\theta)=(r\cos(\theta+\varphi(r)),r\sin(\theta+\varphi(r)))\quad(r>\sqrt{2}),
$$
where $\varphi(r)=\arcsin(2/r^2)$, then for any point $p$ with $d(O,p)>\sqrt{2}$, the points $O, p, f(p)$ form an isosceles triangle of area $1$.
Therefore, it suffices to find an appropriate point $p\in S$ such that $f(p)$ is also contained in $S$. This can be achieved using Lebesgue's density theorem.
For the case of right-angled triangles, we can use the same argument, with $f(p)$ replaced by $g(p)=(p+f(p))/2$.

Noting that isosceles trapezoids are obtained by cutting the apex of an isosceles triangle, a slight modification of the above argument proves the following result:

\begin{theorem}\label{thm:trapezoid}
    Let $S\subset\mathbb{R}^2$ be a measurable set of infinite Lebesgue measure. Then, $S$ contains the four vertices of an isosceles trapezoid of area $1$.
\end{theorem}

We note that \Cref{thm:trapezoid} is not true for general unbounded sets of positive Lebesgue measure.
For example, the set $\{(x,y)\mid x^2+y^2\leq 1/100\}\cup \{(n,0)\mid n=1,2,3,\dots\}$ gives a counterexample.

\subsection*{Notation}
For two points $x,y\in \mathbb{R}^2$, we write $d(x,y)$ for the Euclidean distance between $x$ and $y$.
We write $B(x,r)$ for the Euclidean open disk centered at $x$ with radius $r$.
For a measurable set $S\subset \mathbb{R}^2$, we write $\mu(S)$ for the two-dimensional Lebesgue measure of $S$.


\section{Proof of Theorem 1}

    Fix a sufficiently large real number $C>0$ (e.g. $C=100$).
    By Lebesgue's density theorem, there exists a point $A\in S$ and a real number $\varepsilon\in (0,1)$ such that
    \begin{align}\label{eq:density}
        \dfrac{\mu(B(A,\varepsilon)\cap S)}{\mu(B(A,\varepsilon))}>\dfrac{9}{10}.
    \end{align}
    Since $S$ is unbounded, there exists a point $O\in S$ such that
    $$
        d:=d(O,A)>\dfrac{C}{\varepsilon}+\varepsilon.
    $$
    We may assume that $O=(0,0)$ and $A=(d,0)$.
    Let $B=B(A,\varepsilon)$, $B'=B(A,\varepsilon/2)$, and $D=B(O,C)$.
    We note that $B$ and $D$ are disjoint.

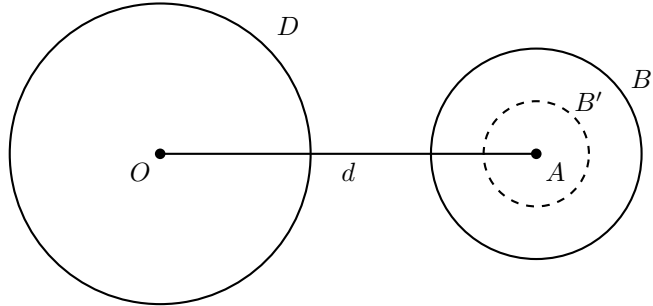
\begin{figure}[h]
\begin{tikzpicture}
  \coordinate (O) at (0,0);
  \coordinate (A) at (5,0);

  \fill[black] (O) circle (2pt) node[below left] {$O$};
  \fill[black] (A) circle (2pt) node[below right] {$A$};

  \draw[thick] (O) circle (2);
  \node at (1.7,1.7) {$D$};

  \draw[thick, dashed] (A) circle (0.7);
  \node at (5.7,0.7) {$B'$};

  \draw[thick] (A) circle (1.4);
  \node at (6.4,1.0) {$B$};

  \draw[thick] (O) -- (A);

  \coordinate (M) at (2.5,0);
  \node[below] at (M) {$d$};

\end{tikzpicture}
\end{figure}

    First we prove the statement for isosceles triangles.
    We define a $C^\infty$-function $\varphi\colon (\sqrt{2},\infty)\to (0,\pi/2)$ by
    $$
    \varphi(r)=\arcsin\left(\dfrac{2}{r^2}\right),
    $$
    and define a $C^\infty$-diffeomorphism $f\colon \mathbb{R}^2\setminus D\to \mathbb{R}^2\setminus D$ by
    $$
        f(r\cos\theta,r\sin\theta)=(r\cos(\theta+\varphi(r)),r\sin(\theta+\varphi(r))).
    $$
    It is clear that the Jacobian determinant of $f$ is given by $J_f(p)=1$ for all $p$.
    Also, for any $p\in \mathbb{R}^2\setminus D$, the triangle formed by the three points $O,p,f(p)$ is an isosceles triangle of area $1$:
    $$
        \dfrac{r^2}{2}\sin(\varphi(r))=\dfrac{r^2}{2}\cdot\dfrac{2}{r^2}=1.
    $$
    Therefore, it suffices to find a point $p\in B'\cap S$ such that $f(p)\in B\cap S$.

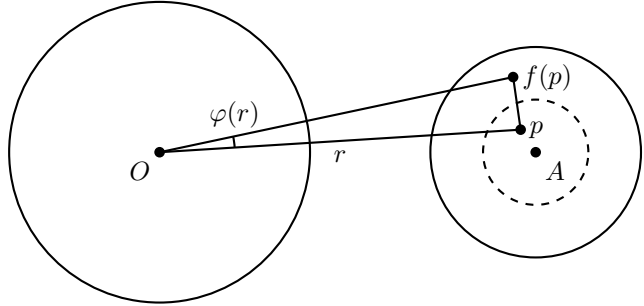
    \begin{figure}[h]
\begin{tikzpicture}
  \coordinate (O) at (0,0);
  \coordinate (A) at (5,0);
  \coordinate (P) at (4.8,0.3);
  \coordinate (M) at (2.4,0.15);
  \coordinate (F) at (4.7,1.0);
  \coordinate (V) at (1.0,0.5);

  \fill[black] (O) circle (2pt) node[below left] {$O$};
  \fill[black] (A) circle (2pt) node[below right] {$A$};
  \fill[black] (P) circle (2pt) node[right] {$p$};
  \fill[black] (F) circle (2pt) node[right] {$f(p)$};
  \node at (V) {$\varphi(r)$};
  \draw[thick] (O) circle (2);
  \node at (1.7,1.7) {};

  \draw[thick, dashed] (A) circle (0.7);
  \node at (5.7,0.7) {};

  \draw[thick] (A) circle (1.4);
  \node at (6.4,1.0) {};

  \draw[thick] (O) -- (P);
  \draw[thick] (F) -- (O);
  \draw[thick] (F) -- (P);

  \node[below] at (M) {$r$};

  \draw pic[draw, thick, angle radius=1cm] {angle=P--O--F};

\end{tikzpicture}
\end{figure}

    Suppose that $p=(r\cos\theta,r\sin\theta)\in B'$.
    By the inequality
    $$
    \sin x<x<\dfrac{\pi}{2}\sin(x)
    $$
    which holds for $x\in (0,\pi/2)$, we have
    $$
    d(p,f(p))=2r\sin\dfrac{\varphi(r)}{2}<r\varphi(r)<\dfrac{r\pi}{2}\sin(\varphi(r))=\dfrac{\pi}{r}.
    $$
    Since we have $r>d-\varepsilon>C/\varepsilon$, we get
    $$
        d(p,f(p))<\dfrac{\pi\varepsilon}{C}<\dfrac{\varepsilon}{2}.
    $$
    Therefore, by the triangle inequality, we have the inclusion $f(B')\subset B$.
    Suppose for contradiction that $f(B'\cap S)\cap S=\varnothing$.
    Then we have $f(B'\cap S)\subset B\setminus S$ and hence
    $$
        \mu(B\setminus S)\geq \mu(B'\cap S)\geq \mu(B')-\mu(B\setminus S).
    $$
    This implies $\mu(B\setminus S)\geq \mu(B')/2=\mu(B)/8$, which contradicts the inequality \eqref{eq:density}.
    Therefore, there is some point $p\in B'\cap S$ such that $f(p)\in B\cap S$.
    This finishes the proof for isosceles triangles.

    To prove the statement for right-angled triangles, it suffices to show that $S$ contains the three vertices of a right-angled triangle of area $1/2$.
    We define an injective $C^\infty$-map $g\colon \mathbb{R}^2\setminus D\to \mathbb{R}^2$ by
    $$
        g(r\cos\theta,r\sin\theta)=\dfrac{p+f(p)}{2}=\left(r\cos\left(\dfrac{\varphi(r)}{2}\right)\cos\left(\theta+\dfrac{\varphi(r)}{2}\right),r\cos\left(\dfrac{\varphi(r)}{2}\right)\sin\left(\theta+\dfrac{\varphi(r)}{2}\right)\right).
    $$
    Then, for any $p\in \mathbb{R}^2\setminus D$, the triangle formed by the three points $O,p,g(p)$ is a right-angled triangle of area $1/2$.
    Therefore, it suffices to find a point $p\in B'\cap S$ such that $g(p)\in B\cap S$.
    Suppose that $p=(r\cos\theta,r\sin\theta)\in B'$.
    Then, we have
    $$
    d(p,g(p))=\dfrac{d(p,f(p))}{2}<\dfrac{\varepsilon}{4}.
    $$
    Therefore, by the triangle inequality, we have the inclusion $g(B')\subset B$.
    The Jacobian determinant of $g$ at $p$ is given by
    $$
        J_g(p)=\cos\left(\dfrac{\varphi(r)}{2}\right)\cdot\dfrac{d}{dr}\left(r\cos\left(\dfrac{\varphi(r)}{2}\right)\right)=\cos^2\left(\dfrac{\varphi(r)}{2}\right)+\dfrac{\sin(\varphi(r))}{\sqrt{r^4-4}}.
    $$
    Since we have $r>C$ and hence $\varphi(r)<\varphi(C)<\pi/4$, we get
    $$
        J_g (p) > \cos^2 \dfrac{\pi}{8}>\dfrac{4}{5}.
    $$
    Suppose for contradiction that $g(B'\cap S)\cap S=\varnothing$.
    Then we have $g(B'\cap S)\subset B\setminus S$ and hence
    $$
        \mu(B\setminus S)\geq \dfrac{4}{5}\mu(B'\cap S)\geq \dfrac{4}{5}(\mu(B')-\mu(B\setminus S)).
    $$
    This implies the inequality
    $$
        \mu(B\setminus S)\geq \dfrac{4}{9}\mu(B')=\dfrac{1}{9}\mu(B),
    $$
    which contradicts the inequality \eqref{eq:density}.
    Therefore, there is some point $p\in B'\cap S$ such that $g(p)\in B\cap S$.
    This finishes the proof for right-angled triangles.

\section{Proof of Theorem 2}

    Fix a sufficiently large real number $C>0$ (e.g. $C=100$).
    By Lebesgue's density theorem, there exists a point $a\in S$ and a real number $\varepsilon\in (0,1)$ such that
    \begin{align*}
        \dfrac{\mu(B(a,\varepsilon)\cap S)}{\mu(B(a,\varepsilon))}>\dfrac{9}{10}.
    \end{align*}
    Since $S$ has infinite Lebesgue measure, again by Lebesgue's density theorem, there exists a point $O\in S$ such that
    \begin{align}\label{eq:density_3}
        \lim_{\delta\to 0^+}\dfrac{\mu(B(O,\delta)\cap S)}{\mu(B(O,\delta))}=1\quad\text{and}\quad d:=d(O,A)>\dfrac{C}{\varepsilon}+\varepsilon.
    \end{align}
    We may assume that $O=(0,0)$ and $A=(d,0)$.
    Let $B=B(A,\varepsilon)$, $B'=B(A,\varepsilon/2)$, and $D=B(O,C)$.
    We note that $B$ and $D$ are disjoint.

    We first replace $S$ with its subset with a nice property.
    For a real number $R>0$, we consider the set
    $$
    S_R:=S\cap (R\cdot S) \subset S,
    $$
    that is, the set of points $p\in S$ such that $R^{-1}p$ is also in $S$.
    Then, we have
    $$
        (B\cap S)\setminus S_R\subset B(O,2d)\setminus (R\cdot S)=R\cdot (B(O,2d/R)\setminus S)
    $$ and hence
    \begin{align*}
    \mu((B\cap S)\setminus S_R)&\leq R^2\cdot\bigl(\mu(B(O,2d/R))-\mu(B(O,2d/R)\cap S)\bigr)\\
    &=R^2\cdot \mu(B(O,2d/R))\cdot \left(1-\dfrac{\mu(B(O,2d/R)\cap S)}{\mu(B(O,2d/R))}\right)\\
    &=\mu(B(O,2d))\cdot \left(1-\dfrac{\mu(B(O,2d/R)\cap S)}{\mu(B(O,2d/R))}\right).
    \end{align*}
    Therefore, by \eqref{eq:density_3}, we have $\lim_{R\to \infty}\mu((B\cap S)\setminus S_R)=0$ and hence $\lim_{R\to \infty}\mu(B\cap S_R)=\mu(B\cap S)$.
    In particular, there exists some $R>100$ such that
    \begin{align}\label{eq:density_4}
        \dfrac{\mu(B\cap S_R)}{\mu(B)}>\dfrac{8}{9}.
    \end{align}
    We fix such $R$ until the end of the proof.

    The rest of the proof is essentially the same as that of \Cref{thm:iso_right} (1).
    We define a $C^\infty$-function $\psi\colon (2,\infty)\to (0,\pi/2)$ by
    $$
    \psi(r)=\arcsin\left(\dfrac{R^2}{R^2-1}\cdot\dfrac{2}{r^2}\right),
    $$
    and define a $C^\infty$-diffeomorphism $f\colon \mathbb{R}^2\setminus D\to \mathbb{R}^2\setminus D$ by
    $$
        f(r\cos\theta,r\sin\theta)=(r\cos(\theta+\psi(r)),r\sin(\theta+\psi(r))).
    $$
    It is clear that the Jacobian determinant of $f$ is given by $J_f(p)=1$ for all $p$.
    Also, for any $p\in \mathbb{R}^2\setminus D$, the quadrilateral formed by the four points $p,f(p),R^{-1}f(p),R^{-1}p$ is an isosceles trapezoid of area $1$:
    $$
        \left(1-\dfrac{1}{R^2}\right)\cdot\dfrac{r^2}{2}\sin(\psi(r))=
        \dfrac{R^2-1}{R^2}\cdot\dfrac{r^2}{2}\cdot\dfrac{R^2}{R^2-1}\cdot\dfrac{2}{r^2}=1.
    $$
    Therefore, it suffices to find a point $p\in B'\cap S_R$ such that $f(p)\in B\cap S_R$.
    Suppose that $p=(r\cos\theta,r\sin\theta)\in B'$.
    As in the proof of \Cref{thm:iso_right}, we have
    $$
        d(p,f(p))<\dfrac{\varepsilon}{2}.
    $$
    Therefore, by the triangle inequality, we have the inclusion $f(B')\subset B$.
    Suppose for contradiction that $f(B'\cap S_R)\cap S_R=\varnothing$.
    Then we have $f(B'\cap S_R)\subset B\setminus S_R$ and hence
    $$
        \mu(B\setminus S_R)\geq \mu(B'\cap S_R)\geq \mu(B')-\mu(B\setminus S_R).
    $$
    This implies $\mu(B\setminus S_R)\geq \mu(B')/2=\mu(B)/8$, which contradicts the inequality \eqref{eq:density_4}.
    Therefore, there is some point $p\in B'\cap S_R$ such that $f(p)\in B\cap S_R$.
    This finishes the proof of \Cref{thm:trapezoid}.

\printbibliography

\end{document}